\documentstyle[12pt,twoside, amscd]{amsart}
\newcommand{\cI}{{\cal I}}

\newcommand{\cH}{{\cal H}}
\newcommand{\cR}{{\cal R}}

\newcommand{\cN}{{\cal N}}

\newcommand{\cE}{{\cal E}}
\newcommand{\cO}{{\cal O}}

\newcommand{\cA}{{\cal A}}

\newcommand{\cD}{{\cal D}}
\newcommand{\e}{{\bf e}}
\newcommand{\Po}{{\bf P}}

\newcommand{\Aut}{\operatorname{Aut}}

\newcommand{\diff}{\operatorname{d}}
\newcommand{\Bl}{\operatorname{Bl}}

\newtheorem{thm}{Theorem}[section]

\newtheorem{lem}[thm]{Lemma}
\newtheorem{cor}[thm]{Corollary}

\theoremstyle{definition}
\newtheorem{defn}[thm]{Definition}

\theoremstyle{remark}

\def\proof{\smallskip{\it Proof. }}
\def\endproof{\hfill\qed}
\def\qed {\nobreak$\quad$\lower 1pt\vbox{\hrule
\hbox to 8pt{\vrule height 8pt\hfil\vrule height 8pt}
\hrule}\ifmmode\relax\else\par\medbreak \fi}

\begin{document}
\title[A Compactification of Open Varieties] 
{
A Compactification of Open Varieties}
\author[Yi {\sc Hu} ]{Yi {\sc Hu}}


\begin{abstract}
In this paper we prove a general method to compactify certain open varieties
by adding normal crossing divisors. This is done by showing that {\it blowing up
 along an  arrangement  of subvarieties} can be carried out. Important examples
such as  Ulyanov's configuration spaces
and complements of arrangements of linear subspaces in projective spaces,  etc.,  are covered.
Intersection ring and (non-recursive) Hodge polynomials
are computed. Furthermore, some general structures
arising from the blowup process are also described and studied.

\end{abstract}
\maketitle




\section{Introduction and the main theorems}

Throughout the paper, the base field is assumed to be algebraically closed.

Let $S$ be a partially ordered set (poset). The rank of $s \in S$ is the maximum of the lengths
of all the chains that end up at $s$. A minimal element is of rank $0$. The rank of $S$ is
the maximum of the lengths of all chains. Let $S_{\le r}$ be the
subposet of elements of rank $\le r$.
All posets in this paper are {\it partially ordered by inclusion} unless otherwise stated.

Two smooth closed subvarieties $U$ and $V$ of a smooth variety $W$ are said to intersect
cleanly if the scheme-theoretic intersection $U \cap V$ is smooth
and $T(U \cap V) = T(U) \cap T(V)$ for their tangent spaces.

\begin{thm}
\label{mainlemma}  Let $X^0$ be an open subset of
a nonsingular algebraic variety $X$. Assume that
$X     \setminus  X^0 $ can be decomposed as a union $ \bigcup_{i \in I} D_i$ of
closed irreducible subvarieties such that
\begin{enumerate}
\item $D_i$ is smooth;
\item $D_i$ and $D_j$ meet cleanly;
\item $D_i \cap D_j = \emptyset$ or a disjoint union of $D_l$.
\end{enumerate}
The set $\cD = \{D_i\}_i$ is then a poset. Let $k$ be the rank of $\cD$.
Then there is a sequence of well-defined blowups
$$\Bl_{\cD} X  \rightarrow \Bl_{\cD_{\le k-1}} X
 \rightarrow \ldots \rightarrow \Bl_{\cD_{\le 0}} X \rightarrow  X$$
  where $\Bl_{\cD_{\le 0}} X \rightarrow  X$ is the  blowup
of $X$ along $D_i$ of rank 0,
and inductively, $\Bl_{\cD_{\le r}} X \rightarrow \Bl_{\cD_{\le r-1}} X$
is the blowup of $\Bl_{\cD_{r-1}} X$ along the proper transforms of $D_j$ of rank $r$, such that
\begin{enumerate}
\item $\Bl_{\cD} X$ is smooth;
\item $\Bl_\cD X     \setminus  X^0 = \bigcup_{i \in I} \widetilde{D}_i$ is a divisor with normal
crossings
\item $\widetilde{D}_{i_1} \cap \ldots \cap \widetilde{D}_{i_n}$ is nonempty
if and only if $D_{i_1} \ldots D_{i_n}$ form a chain in the poset $\cD$.
Consequently,  $\{\widetilde{D}_i\}$ and $\{\widetilde{D}_j\}$  meet
if and only if   $\{D_i\}$ and $\{D_j\}$ are comparable.
\end{enumerate}
\end{thm}

\begin{defn}
\label{arr} The set $\cD$
is called an {\it arrangement} of smooth subvarieties.
$\Bl_\cD X$ is referred as the blowup of $X$ {\it along the arrangement} $\cD$ of subvarieties.
When the condition (3) of  Theorem \ref{mainlemma}  is replaced by:
($3'$). $D_i \cap D_j = \emptyset$ or $D_l$ for some $l$,
then $\cD$ is called  a {\it simple arrangement} of smooth subvarieties.
\end{defn}

This work naturally extends the previous works of Fulton-MacPherson (\cite{FM}),
MacPherson-Procesi \cite{MP} and Ulyanov (\cite{Ul}). Our main
theorem was especially inspired by Ulyanov's paper (\cite{Ul}).

\medskip

Any collection of (affine) linear subspaces $\{H_i\}$
in ${\Bbb P}^n$ (or ${\Bbb C}^n$) induces a simple arrangement of smooth subvarieties
by taking all possible non-empty intersections (subspace arrangement).
  Theorem \ref{mainlemma} applies to such situation.
A smooth curve of higher degree and a general line  in ${\Bbb P}^2 \subset {\Bbb P}^n$
($n>2$) necessarily meet in several distinct points.
Hence it is useful to include as well  non-simple arrangements of subvarieties.

More  sophisticated and important examples are in order
to situate Theorem \ref{mainlemma} in particular cases followed by stating
certain general structures arising from the construction
of the blowup $\Bl_\cD X$.

\smallskip\noindent
{\sl  $\bullet$ Configuration spaces.}
\smallskip

Consider $X^n$. Let $\Delta_{ij}$ be the subset of all points whose $i$-th and $j$-th
coordinates coincide. Let $\Delta$ be the set of all possible intersections of $\Delta_{ij}$.
$\Delta$ satisfies the arrangement conditions.
We will call $\Delta$ the diagonal arrangement.

\begin{cor} (Ulyanov \cite{Ul})
 $\Bl_\Delta (X^n)$ is a  symmetric\footnote{i.e., $\Sigma_n = \Aut \{1, \ldots, n\}$
acts on it} smooth
projective compactification of $X^n     \setminus  \bigcup \Delta$ by adding
smooth divisors with normal crossings.
\end{cor}

Ulyanov also proved that  $X\langle n \rangle := \Bl_\Delta (X^n)$
dominates $X [n]$,
the Fulton-MacPherson configuration space (\cite{FM}).
$X[n]$ is not an instance of blowups along {\it arrangement} of subvarieties
(cf. Definition \ref{arr}).

\smallskip\noindent
{\sl  $\bullet$ Space of holomorphic maps.}
\smallskip

Let $N_d ({\Bbb P}^n)$ be the space of  $(n+1)$-tuples $(f_0, \ldots, f_n)$ modulo homothety
where $f_i$ are homogeneous polynomials of degree $d$ in two variables and
 $N_d^0({\Bbb P}^n)$  the open subset such that $f_0, \ldots, f_n$ have no common zeros.
 $N^0_d({\Bbb P}^n)$ is
 naturally identified with  the  space $M^0_d({\Bbb P}^n)$ of holomorphic maps of degree $d$
from ${\Bbb P}^1$ to ${\Bbb P}^n$.

For any  integer $0 \le d_0 \le d$ and  an arbitrary partition $\tau= \sum_j d_j$ of $d-d_0$,
let $$N_{d_0, \tau} :=\{ [\sigma_0 \prod_{j >0} (a_j w_0 - b_j w_1)^{d_j}
 \ldots \sigma_n \prod_{j >0} (a_j w_0 - b_j w_1)^{d_j}] \} \subset N_d$$
where  $\sigma_j$ are homogeneous polynomials of degree $d_0$
and $\{[a_j, b_j]\}_j$ are unordered points in ${\Bbb P}^1$.

The collection $\cN = \{N_{\pi, d_0}\}$ is however not an arrangement of smooth subvarieties in $N_d$,
thanks to an important observation by Sean Keel who saves me an embarrassment.
For, some strata $N_{\pi, d_0}$ may have singularities along  lower strata. But the minimal ones are smooth
so that the first step of the iterated blowups can be carried out. The hope is that after the first step,
the singularieties of the strata on the second level are resolved and get separated so that the next
step of the iterated blowups can also be carried out. It calls for further investigation
 to see if the process can indeed be 
excuted step by step to obtain a good compactification $\Bl_{\cN} N_d({\Bbb P}^n)$ 
of the space $N_d^0({\Bbb P}^n)$   of holomorphic maps.


\smallskip\noindent
{\sl  $\bullet$ GIT.}
\smallskip

Theorem \ref{mainlemma} coupled with a compatible group action yields
an instance of Theorem \ref{mainlemma} in Geometric Invariant Theory.
Roughly, it says that blowing up along an arrangement descends to blowing up
of any GIT quotient along an induced arrangement.
As a particular case, we recover
Kirwan's  partial desingularization of singular  GIT quotients.
 See Corollary \ref{useform0nbar} in \S \ref{GIT}.

\medskip

We now return to the general situation.

\smallskip\noindent
{\sl  $\bullet$  Proper transforms and exceptional divisors.}
\smallskip

Of useful computational value is that
in each stage of the blowups,  $\Bl_{\cD_{\le r}} X \rightarrow X$,
 the proper transforms of $D_i$ and exceptional divisors are special instances
of Theorem \ref{mainlemma} and all are concisely described using posets induced from $\cD$.

\begin{thm}
\label{transforms} (Proper transforms)
Let $\cD^{r+1}$ be the set of proper transforms in $\Bl_{\cD_{\le r}} X$ of $D_i$ of rank $\ge r+1$. Then
\begin{enumerate}
\item
$$\cD^{r+1} = \{ \Bl_{(\cD_{< D_i})_{\le r}} D_i : \hbox{rank} (D_i) \ge r+1 \}$$
where $\cD_{< D_i}$ denotes the subposet of the elements less than $D_i$, and
$$\cD^{r+1}_{\le 0} = \{ \Bl_{(\cD_{< D_i})_{\le r}} D_i : \hbox{rank} (D_i) = r+1 \}$$
\item $\cD^{r+1}$ is an arrangement of smooth subvarieties in $\Bl_{\cD_{\le r}} X$.
\end{enumerate}
\end{thm}

\begin{cor}
\label{iterated}
With the above notations, $\Bl_\cD X $ can be expressed as iterated blowups along
     (explicit) disjoint centers.
$$\Bl_\cD X = \Bl_{\cD^k_{\le 0}} \Bl_{\cD^{k-1}_{\le 0}} \ldots \Bl_{\cD^1_{\le 0}} \Bl_{\cD_{\le 0}} X.$$
\end{cor}

\begin{cor}
\label{intermediatestep}
The intermediate blowup
$$\Bl_{\cD_{\le r}} X \rightarrow X$$ is  an instance of the theorem
when the arrangement of subvarieties is the sub-arrangement $\cD_{\le r}$.
In particular,
$$\Bl_{\cD_{\le r}} X = \Bl_{\cD^r_{\le 0}} \Bl_{\cD^{r-1}_{\le 0}} \ldots \Bl_{\cD^1_{\le 0}} \Bl_{\cD_{\le 0}} X.$$
\end{cor}


\begin{thm}
\label{divisors}
(Exceptional divisors) Let $\cE^{r+1}$ be the set of all the
 exceptional divisors of $\Bl_{\cD_{\le r}} X \rightarrow X$. Then $\cE^{r+1}$ consists of
\begin{enumerate}
\item for each $D_i$ of rank $r$, $$E^{r+1}_i = {\Bbb P}(N_{D_i^r/\Bl_{\cD_{\le r-1} X}}),$$
 These are also exceptional divisors of the blowup
$\Bl_{\cD_{\le r}} X \rightarrow \Bl_{\cD_{\le r-1}} X$, where $\Bl_{\cD_{\le -1}} X := X$;
\item for each $D_i$ of rank $m < r$,
$$E^{r+1}_i = \Bl_{\{{D^{m+1}_j \cap E_i^{m+1} : D_j > D_i}\}_{\le r-m}} E^{m+1}_i.$$
Note that by (1), $E^{m+1}_i = {\Bbb P} (N_{D_i^m/\Bl_{\cD_{\le m-1}} X}).$
\end{enumerate}
\end{thm}

The expression of $E^{r+1}_i$ relies on the proper transform $D_i^m$ and
is the blowup of ${\Bbb P} (N_{D_i^m/\Bl_{\cD_{\le m-1}} X})$ along an induced arrangement of
smooth subvarieties. Some topological calculations on $N_{D_i^m/\Bl_{\cD_{\le m-1}} X}$
can be reduced to $N_{D_l/X}$ for $D_l \le D_i$ (e.g., \S \ref{rings}).

\smallskip\noindent
{\sl  $\bullet$ Intersection rings.}
\smallskip

An embeding $U \hookrightarrow W$ is called a Lefschetz embedding if the restriction map
$A^\bullet (W) \rightarrow A^\bullet (U)$ is surjective.
Let $J_{U/W}$ be the kernel of $A^\bullet (W) \rightarrow A^\bullet (U)$
and $P_{U/W}$ be a Chern polynomial for the normal bundle $N_{U/W}$ (\S \ref{rings}).

A simple arrangement $\cD$ is called {\it regular} if for any $D_l < D_i$ there is
$D_j > D_l$ such that $D_l = D_i \cap D_j$. That is, any $D_l$ is an intersection of
maximal $D_i$.  All the previously mentioned examples
are regular. 

\begin{thm}
\label{ring}
Let $\cD$ be a regular simple arrangement of subvarieties.
  Assume that all inclusion $D_i \subset D_j$ and $D_i \subset X$ are Lefschetz embeddings.
Then the Chow ring $A^\bullet (\Bl_\cD X)$ is isomorphic to the  polynomial ring
$$A^\bullet (X) [T_1, \ldots T_N]/I$$
where $T_i$ corresponds to\footnote{
More precisely and geometrically,
$T_i$ corresponds to the exceptional divisor $E^{r+1}_i$ where, $r= \hbox{rank} (D_i)$. See
the proof in \S \ref{rings}.
See also Theorem \ref{divisors} (1) for the description of $E^{r+1}_i$.}
 $D_i$
and $I$ is the ideal generated by
\begin{enumerate}
\item $T_i \cdot T_j$, for incomparable $D_i$ and $D_j$;
\item $J_{D_i/X} \cdot T_i$; for all $i$;
\item $P_{D_i/X} (- \sum_{D_j \le D_i} T_j)$ for all $i$.
\end{enumerate}
\end{thm}

This theorem can be directly applied to the intermediate stage $\Bl_{\cD_{\le r}} X$ by Corollary
\ref{intermediatestep}, to the proper transforms $D^r_i$ by Theorem \ref{transforms},
to the exceptional divisors by Theorem \ref{divisors}.

\smallskip\noindent
{\sl   $\bullet$  Hodge and Poincar\'e polynomials.}
\smallskip

The concise presentations of $\Bl_\cD X$, the intermediate stage $\Bl_{\cD_{\le r}} X$
and the proper transform $D^r_i$ of $D_i$ in every  stage allow one to derive a concise
{\it non-recursive} formula  for the Hodge (Poincar\'e) polynomial of $\Bl_\cD X$.

Let $\e (W)$ ($\Po (W)$) be the Hodge (Poincar\'e) polynomial in two (resp. one)
variables $u$ and $v$ (resp. $t$)
of a smooth projective variety $W$.

\begin{thm}
\label{hodge}
$$\e (\Bl_\cD X) = \e (X) + \sum_{{\scriptsize \begin{array}{c}
D_{i_1} < \ldots <D_{i_{r+1}}\cr
 D_{i_{r+1}} := X
\end{array}}}
\e (D_{i_1}) \prod_{j=1}^r {{(uv)^{\dim D_{i_{j+1}} - \dim D_{i_j}}} - uv \over {uv -1}}.$$
$$\Po (\Bl_\cD X) = \Po (X) + \sum_{{\scriptsize \begin{array}{c}
D_{i_1} < \ldots < D_{i_{r+1}} \cr
 D_{i_{r+1}} := X \end{array} }}
\Po (D_{i_1}) \prod_{j=1}^r {{t^{2\dim D_{i_{j+1}} - 2\dim D_{i_j}}} - t^2 \over {t^2 -1}}.$$
\end{thm}

This formula immediately  applies to the intermediate stage $\Bl_{\cD_{\le r}} X$,
the proper transform $D^r_i$ and the exceptional divisor $E^r_i$.

Consider $X\langle n \rangle$.
The index set of the diagonal arrangement $\Delta$ is the set
 of all partitions $\pi$ of $[n] = \{ 1, \ldots n\}$ except the
largest trivial partition
${\bf 1}_{[n]} = 1 \cup \ldots \cup n.$
The subvariety $\Delta_\pi \in \Delta$ is the set of all points any two of whose
coordinates coincide whenever their indexes belong to the same block of $\pi$.
 Let $\rho (\pi)$ be the number of blocks of $\pi$. Then $\Delta_\pi \cong X^{\rho (\pi)}$.

\begin{cor}
\label{config}
$$\e (X \langle n \rangle) = \e (X)^n + \sum_{{\scriptsize \begin{array}{c}
\pi_{i_1} < \ldots <\pi_{i_{r+1}} \cr
 \pi_{i_{r+1}} : ={\bf 1}_{[n]} \end{array}}}
\e (X)^{\rho (\pi_{i_1})} \prod_{j=1}^r
{{(uv)^{\dim X ( \rho (\pi_{i_{j+1}})  -  \rho(\pi_{i_j}))}} - uv \over {uv -1}}.$$
$$\Po (X \langle n \rangle) = \Po (X)^n + \sum_{{\scriptsize \begin{array}{c}
\pi_{i_1} < \ldots <\pi_{i_{r+1}} \cr
 \pi_{i_{r+1}} : ={\bf 1}_{[n]} \end{array}}}
\Po (X)^{\rho (\pi_{i_1})} \prod_{j=1}^r
 {{(t^2)^{\dim X (\rho (\pi_{i_{j+1}})  - \rho(\pi_{i_j}))}} - t^2 \over {t^2 -1}}.$$
\end{cor}

\begin{cor}
\label{hyperplane}
 Let $\cH = \{ H_i\}$ be an arrangement of linear subspaces of ${\Bbb P}^n$. Then
 $$\e (\Bl_\cH {\Bbb P}^n) = {1 \over uv}
\sum_{{\scriptsize \begin{array}{c}
 H_{i_0} < H_{i_1} < \ldots <H_{i_{r+1}} \cr
 \emptyset := H_{i_0},         H_{i_{r+1}} := X \end{array}}}
 \prod_{j=0}^r {{(uv)^{\dim H_{i_{j+1}} - \dim H_{i_j}}} - uv \over {uv -1}},$$
$$\Po (\Bl_\cH {\Bbb P}^n) = {1 \over t^2}
\sum_{{\scriptsize \begin{array}{c}
 H_{i_0} < H_{i_1} < \ldots <H_{i_{r+1}} \cr
\emptyset := H_{i_0},  H_{i_{r+1}} := X \end{array}}}
 \prod_{j=0}^r {{t^{2\dim H_{i_{j+1}} - 2\dim H_{i_j}}} - t^2 \over {t^2 -1}},$$
where $\dim \emptyset = - 2$.
\end{cor}

Take $n+2$ points of ${\Bbb P}^n$ in general linear position.
 They span $n+2 \choose 2$ hyperplanes. Let $\cH_{n}$ be the induced simple arrangement.
Then $\Bl_{\cH_n} {\Bbb P}^n$ is isomorphic to $\overline{M}_{0,n+3}$. This example is due to Kapranov.
The index set of $\cH_{n}$ is the set of all subsets $S$ of $[n+2]$ such that $1 \le |S| \le n$.

\begin{cor}
\label{m0nhodge}
 $$\e (\overline{M}_{0,n+3}) = {1 \over uv}
\sum_{{\scriptsize  \begin{array}{c}
 S_{i_0} < S_{i_1} < \ldots <S_{i_r} < S_{i_{r+1}} \cr
\emptyset := S_{i_0} ,|S_{i_{r+1}}| := n+1 \end{array}} }
 \prod_{j=0}^r {{(uv)^{ |S_{i_{j+1}}| - |S_{i_j}|}} - uv \over {uv -1}},$$
$$\Po (\overline{M}_{0,n+3}) = {1 \over t^2}
\sum_{{\scriptsize  \begin{array}{c}
 S_{i_0} < S_{i_1} < \ldots <S_{i_r} < S_{i_{r+1}} \cr
\emptyset := S_{i_0} ,|S_{i_{r+1}}| := n+1 \end{array}} }
 \prod_{j=0}^r {{t^{2|S_{i_{j+1}}| - 2| S_{i_j}}|} - t^2 \over {t^2 -1}},$$
where $|\emptyset| = - 1$.
\end{cor}

Keel computed these numbers and furthermore he also computed the intersection ring
( \cite{keel}).

\begin{cor}
\label{inCn}
 Let $\cH = \{ H_i\}$ be an arrangement of linear subspaces of ${\Bbb C}^n$. Then
$$\e (\Bl_\cH {\Bbb C}^n) = 1 + {1 \over uv}
\sum_{{\scriptsize  \begin{array}{c}
H_{i_1} < \ldots < H_{i_r} < H_{i_{r+1}} \cr
 H_{i_{r+1}} := {\Bbb C}^n \end{array}}
}
 \prod_{j=1}^r {{(uv)^{\dim H_{i_{j+1}} - \dim H_{i_j}}} - uv \over {uv -1}}.$$
$$\Po (\Bl_\cH {\Bbb C}^n) = 1 + {1 \over t^2}
\sum_{{\scriptsize  \begin{array}{c}
H_{i_1} < \ldots < H_{i_r} < H_{i_{r+1}} \cr
 H_{i_{r+1}} := {\Bbb C}^n \end{array}}
}
 \prod_{j=1}^r {{t^{2\dim H_{i_{j+1}} - 2\dim H_{i_j}}} - t^2 \over {t^2 -1}}.$$
\end{cor}


Finally, it needs to be pointed out that the general procedure of blowing up along
arrangements can be extended to some singular cases as well. This may be necessary
in certain applications (see \S\S 6 and 7).

The paper is structured as follows. \S 2  provides proofs of the statements in this introduction
 on  the {\it structures} of the blowup along arrangement $\Bl_\cD X$,
the exceptional divisors $E^r_i$ and proper transforms $D^r_i$ of $D_i$. Some corollaries
to the proofs are also drawn.
\S 3 gives an alternative construction of
$\Bl_\cD X$ as the closure of the open subset.
 \S 4 proves the statement on the intersection ring of $\Bl_\cD X$.
 \S 5 proves the formulas for Hodge and Poincar\'e polynomials as stated in this introduction.
\S 6 is devoted to the  spaces of  holomorphic maps.
\S  7  treats blowups of  GIT quotients along induced arrangements.

I learnt from Professor Fulton that  Dylan Thurston  (while still an undergraduate
at Harvard) noticed several years
ago $X[n]$ could be constructed by a sequence of symmetric blowups -- but one has to blow up along
ideal sheaves.  The point is
that one can blow up along two smooth subvarieties that meet excessively in a smooth subvariety
without first blowing up the small variety.

I  wonder if $X[n]$ is the {\it minimal symmetric} compactification of
the configuration space $X^n \setminus  \bigcup \Delta$ by adding normal crossing
divisors. I thank Fulton and MacPherson for their powerful original inspiring work \cite{FM}.
 This paper is dedicated to them.

\medskip

{\sl Acknowledgements}.
Our paper clearly follows the ideas and methods of some earlier
works, especially those of Ulyanov \cite{Ul} and MacPherson and Procesi \cite{MP}.
The computation of the intersection ring follows that of Fulton and MacPherson \cite{FM}.
I thank them all. I am very grateful to Professor Fulton for his instructive
comments and generous advice.
 I thank MPI in Bonn for financial support (summer 1999) while this paper was being written.
I am very indebted to Professor S.-T Yau for his valuable support and Professor S. Keel
for pointing out a serious mistake. The research is partially
supported by NSF and NSA.

\section{Proof of Theorems \ref{mainlemma},  \ref{transforms}, \ref{divisors}}

\begin{lem}
\label{basic} Let $U$ and $V$ be two smooth closed subvarieties of a smooth
variety $W$ that intersect cleanly. Then
\begin{enumerate}
\item the proper transforms of $U$ and $V$ in $\Bl_{U \cap V} W$ are disjoint;
\item the proper transform of $V$ in $\Bl_U W$ is isomorphic to $\Bl_{U \cap V} V$;
\item if $Z$ is a smooth subvariety of $U \cap V$, then the proper transforms of
$U$ and $V$ in $\Bl_Z W$ intersect cleanly.
\end{enumerate}
\end{lem}
\proof
All follow from standard arguments.
\endproof


\begin{lem} (Flag Blowup Lemma. \cite{FM} and \cite{Ul})
Let $V^1_0 \subset V^2_0 \subset \ldots \subset V^s_0 \subset W$ be a flag of smooth
subvarieties in a smooth algebraic variety $W_0$. For $k = 1, \ldots, s$,
define inductively: $W_k$ is the blowup of $W_{k-1}$ along $V^k_{k-1}$;
$V^k_k$ is the exceptional divisor in $W_k$; and $V^i_k, k \ne i$, is the proper transform
of $V^i_{k-1}$ in $W_k$. Then the preimage of $V^s_0$ in the resulting variety
$W_s$ is a normal crossing divisors $V_s^1 \cup \ldots V^s_s$.
\end{lem}
\proof See \cite{Ul}.
\endproof

\smallskip\noindent
{\sl Proof of Theorem \ref{mainlemma}.}

\proof Without the awkward but routine verification of the inductive proof,
 the construction goes quite transparently
by the clarification as follows.

First, $\Bl_{\cD_{\le 0}} X \rightarrow X$ is the
blow up of $X$ along the disjoint smooth subvarieties of $D_i$ of rank 0.

Let $D_j^1$ be the proper transform of $D_j$ of rank $\ge 1$.
By Lemma \ref{basic} (1), the proper transforms $D^1_j$ of $D_j$ of rank 1 are disjoint in
$\Bl_{\cD_{\le 0}} X$.
By Lemma \ref{basic} (2) and (3), all $D_j^1$ are smooth and intersect cleanly (or trivially).
If $D_i \cap D_j= \coprod_l D_l$, then $ D_i^1 \cap D^1_j = \coprod_{\hbox{rank} (D_l) > 0} D^1_l$. Otherwise,
$D_i^1 \cap D^1_j = \emptyset$. This shows that
$$\cD^1 = \{ D^1_j = \Bl_{\cD_{<D_j}} D_j : \hbox{rank} (D_j) \ge 1\}$$
is an arrangement of subvarieties in $\Bl_{\cD_{\le 0}} X$.
Moreover,
$$\cD^1_{\le 0} = \{ D^1_j = \Bl_{\cD_{<D_j}} D_j : \hbox{rank} (D_j) = 1\}$$
This grants the next step possible,
which is essentially a repetition of the first step:
$$\Bl_{\cD_{\le 1}} X = \Bl_{\cD^1_{\le 0}} (\Bl_{\cD_{\le 0}} X ) \rightarrow \Bl_{\cD_{\le 0}} X .$$
Note that  rank ($\cD^1$) = rank ($\cD$) $-1$.

Let $\cD^2$ be the proper transform of $D_j$ in $\Bl_{\cD_{\le 1}} X$ of rank $\ge 2$.
The same reasoning as above shows
that  $$\cD^2 = \{ D^1_j = \Bl_{({\cD_{<D_j}})_{\le 1}} D_j : \hbox{rank} (D_j) \ge 2\}$$
is an arrangement of smooth subvarieties and
$$\cD^2_{\le 0} = \{ D^1_j = \Bl_{({\cD_{<D_j}})_{\le 1}} D_j : \hbox{rank} (D_j) = 2\}.$$
Blowing up subvarieties in $\cD^2_{\le 0}$, we obtain
$$\Bl_{\cD_{\le 2}} X = \Bl_{\cD^2_{\le 0}} (\Bl_{\cD_{\le 1}} X ) \rightarrow \Bl_{\cD_{\le 1}} X .$$
Note that  rank ($\cD^2$) = rank ($\cD^1$) $-1$ = rank ($\cD$) $- 2$.

The above can be repeated until the subvarieties in the rank 0 poset  $\cD^k$ are blown up.
That is, the resulting variety from the last step is the iterated blowup along smooth disjoint centers
$$\Bl_\cD X = \Bl_{\cD^k_{\le 0}} \Bl_{\cD^{k-1}_{\le 0}} \ldots \Bl_{\cD^1_{\le 0}} \Bl_{\cD_{\le 0}} X.$$

Statement (1) follows from this description.

If $D_i \cap D_j \ne \emptyset, D_i , D_j$,
that is, $D_i$ and $D_j$ are incomparable,
 by Lemma \ref{basic} (1),
their proper transforms become disjoint at the stage $$\Bl_{\cD_{\le r}} X \rightarrow X$$
for $r = \max \{ \hbox{rank} (D_l) : D_l \subset D_i \cap D_j \}$. Hence
$\widetilde{D}_{i_1} \cap \ldots \cap \widetilde{D}_{i_n}$ is nonempty
if and only if $D_{i_1}, \ldots,  D_{i_n}$ form a chain in the poset $\cD$.
This proves the statement in (3)

Statement (2) then  follows directly from the Flag Blowup Lemma.
Here one needs to observe that for any maximal chain $D_{i_1} < \ldots <D_{i_n}$,  by the above proof of (3),
blowing up the proper transform of any $D_j$ which is
not in the chain
is irrelevant to the intersection $\widetilde{D}_{i_1} \cap \ldots \cap \widetilde{D}_{i_n}$.
Hence the Flag Blowup Lemma applies.
\endproof

We now draw an easy consequence.
Let $\gamma$ be a chain
$$D_{i_1} < \ldots < D_{i_n}$$ and $S_\gamma$ the
intersection $\widetilde{D}_{i_1} \cap \ldots \cap \widetilde{D}_{i_n}$.
Set $$S^0_\gamma := S_\gamma     \setminus  \bigcup_{\gamma' \supset \gamma} S_{\gamma'}.$$
We allow $\gamma = \emptyset$ and define $S_\emptyset := X$ and hence $S^0_\emptyset = X^0$.
Then the normal crossing property implies that

\begin{cor}
$\bigcup_{\gamma} S^0_\gamma$ is a Whitney stratification of $\Bl_{\cD} X$ by
locally closed smooth subvarieties.
\end{cor}

\smallskip\noindent
{\sl Proof of Theorem \ref{transforms}}

\proof
(1) and (2) will be proved simultaneously by using induction on $r$.

When $r=0$ (the case of $\cD^1$), the proof is contained in the proof of Theorem 1.3.

Assume that Statements (1) and (2) are valid for $\cD^r$.

Consider the blowup $\Bl_{\cD^r_{\le 0}} (\Bl_{\cD_{\le r-1}} X) \rightarrow \Bl_{\cD_{\le r-1}} X$.
$\cD^{r+1}$ is  the set of proper transforms $D^{r+1}_i$
 of $D^r_i = \Bl_{(\cD_{<D_i})_{\le r-1}} D_i \in \cD^r$
for $D_i$ of rank $\ge r+1$. Hence
$$D^{r+1}_i = \Bl_{\{D^r_j: D_j < D_i, \;\hbox{rank} (D_j) = r\}}
 \Bl_{(\cD_{<D_i})_{\le r-1}} D_i = \Bl_{(\cD_{<D_i})_{\le r}}D_i.$$

The  reasoning in the proof of Theorem 1.3 for that $\cD^1$ is an arrangement of smooth subvarieties
can be applied to the blowup
$$\Bl_{\cD^r_{\le 0}} (\Bl_{\cD_{\le r-1}} X) \rightarrow \Bl_{\cD_{\le r-1}} X$$
to yield the same statement for $\cD^{r+1}$.

This completes the proof.
\endproof

Specializing to the diagonal arrangement of $X^n$, we draw a sample consequence.
Let $\pi$ be a non-trivial partition of $[n]$. Then we have
\begin{cor}
The proper transform $\Delta_\pi^{\rho (\pi) -1} \cong X \langle \rho (\pi) \rangle$.
\end{cor}

\smallskip\noindent
{\sl Proofs of Corollaries  \ref{iterated} and \ref{intermediatestep}.}

\proof
Corollary \ref{iterated} is contained in the proof of Theorem \ref{mainlemma}.
(This corollary does not logically depend on Theorem \ref{transforms} but
depends on the notations introduced there. To keep the introduction coherent,
we put the statement after Theorem \ref{transforms}.)

Corollary  \ref{intermediatestep} follows from essentially the same reason.
\endproof

\smallskip\noindent
{\sl Proof of Theorem \ref{divisors}.}

\proof
(1) and (2) will be proved simultaneously by using induction on $r$.

When $r=0$ (the case of $\cE^1$),
consider the blowup $\Bl_{\cD_{\le 0}} X \rightarrow X$, the statements are standard.

Assume that the statement is valid for $\cE^r$.

Consider the blowup $\Bl_{\cD^r_{\le 0}} (\Bl_{\cD_{\le r-1}} X) \rightarrow \Bl_{\cD_{\le r-1}} X$.
The center of the blowup are $D^r_i$ for $D_i$ of rank $r$. Hence statement (1) is standard.

The rest of the exceptional divisors of $\Bl_{\cD_{\le r}} X \rightarrow X$ come from
the proper transforms of $E^r_i$ for $D_i$ of rank $m \le  r-1$.
Hence they are
$$E^{r+1}_i = \Bl_{\{E^r_i \cap D^r_l: D_l > D_i \; \hbox{rank} (D_l) = r\}} E^r_i = $$
$$\Bl_{\{E^r_i \cap D^r_l : D_l > D_i \; \hbox{rank} (D_l) = r \}}
(\Bl_{\{D_j^{m+1} \cap E^{m+1}_i: D_j > D_i\}_{\le r-1-m}} E^{m+1}_i).$$
Now observe that
$$\{D_j^{m+1} \cap E^{m+1}_i: D_j > D_i \}^{r-m}_{\le 0} =
\{E^r_i \cap D^r_l :\; D_l > D_i, \; \hbox{rank} (D_l) = r \}.$$
Hence, by Theorem 1.1 (or its proof),
$$E^{r+1}_i = \Bl_{\{D_j^{m+1} \cap E^{m+1}_i: D_j > D_i\}_{\le r-m}} E^{m+1}_i.$$
\endproof

\section{$\Bl_{\cD} X$ as a closure}

\begin{thm}
 $\Bl_{\cD} X$ is the closure of $X^0 = X     \setminus  \bigcup_i D_i$
in $$X \times \prod_i \Bl_{D_i} X.$$
\end{thm}
\proof We will prove the following statement by induction.

$X_{r+1} = \Bl_{\cD_{\le r}} X$ is the closure of $X^0$
in $X \times \prod_{D_i \in \cD_{\le r}} \Bl_{D_i} X$.

When $r=0$, the statement is clear because $X_1$ is the blowup of $X$ along disjoint
$D_i \in \cD_{\le 0}$ and is thus the same as the closure of the graph of the
rational map $X \rightarrow \prod_{D_i \in \cD_{\le 0}} \Bl_{D_i} X$.
Assume that the statement for $X_r$ is proved.
$X_{r+1}$ is the blowup of $X_r$ along the minimal subvarieties in $\cD^r_{\le 0}$.
This is interpretated as blowing up the ideal of sheaf
$$\prod_{\hbox{rank} D_i = r} \cI (D_i^r),$$
where $\cI (D)$ is the ideal sheaf for a closed subvariety $D$.
Denote by $p_r$ the projection $X_r \rightarrow X$ and by $(p_r^* \cI (D_j))$
the ideal generated by the pull-back $p_r^* \cI (D_j)$.
Then by Corollary \ref{intermediatestep}
and Theorem \ref{mainlemma},
$$(p_r^* \cI (D_i)) =\cI (D^r_i) \prod_{D_j < D_i} \cI (E_j^r).$$
Observe that $\prod_{D_j < D_i} \cI (E_j^r)$ is an invertible sheaf.
Hence $$\prod_{\hbox{rank} D_i = r} \cI (D_i^r) \;\; \hbox{and} \;\;
\prod_{\hbox{rank} D_i = r} (p_r^* \cI (D_j)$$
are differed by a multiple of an invertible sheaf.
Therefore blowing up $$\prod_{\hbox{rank} D_i = r} \cI (D_i^r)$$ is
the same as  taking the closure of the graph of
the rational map $X_r \rightarrow \prod_{\hbox{rank} D_i = r} \Bl_{D_i} X$.
By inductive assumption, it is the same as the closure of $X^0$ in
$$X \times \prod_{D_i \in \cD_{\le r-1}} \Bl_{D_i} (X)
\times \prod_{\hbox{rank} (D_i) = r} \Bl_{D_i} X.$$
This finishes the inductive proof.

The statement of the theorem is the case when $r =$ rank ($\cD$).
\endproof

\section{Intersection ring of $\Bl_\cD X$}
\label{rings}

For any inclusion $U \hookrightarrow W$ of a smooth closed subvariety $U$ in a smooth variety $W$,
$J_{U/W}$ denotes the kernel of $$A^\bullet (W) \rightarrow A^\bullet (U).$$
Assume that  $U \hookrightarrow W$ is a Lefschetz embedding, that is,
$A^\bullet (W) \rightarrow A^\bullet (U)$ is surjective. Then $A^\bullet (U) = A^\bullet (W)/ J_{U/W}$.
Define a Chern polynomial $P_{U/W} (t)$ to be a polynomial
$$P_{U/W} (t) = t^d + a_1 t^{d-1} + \cdots + a_{d-1} t + a_d \in A^\bullet W [t],$$
where $d$ is the codimension of $U$ in $W$ and $a_i \in A^i (W)$ is a class whose restriction in $A^i U$
is the Chern class $c_i (N_{U/W})$, where $N_{U/W}$ is the normal bundle of $U$ in $W$.
In addition, it is required that $a_d = [U]$ be the class of $U$, which is a class restricting to
the top Chern class $c_d (N_{U/W})$.

\begin{lem} (\cite{keel})
\label{keel}
Let $\{U_i\}$ be disjoint smooth closed subvarieties of a smooth variety $W$.
Assume that all inclusions  $U_i \hookrightarrow W$ are Lefschetz embeddings.
Then the Chow ring $A^\bullet (\Bl_{\{U_i\}} W)$ is isomorphic to
$$A^\bullet (W) [ T_1, \ldots, T_m]/I$$
where $T_i$ corresponds to  the exceptional divisor $\widetilde{U_i}$ for $U_i$ and
$I$ is the ideal generated by
\begin{enumerate}
\item $T_i \cdot T_j$, for $ i \ne j$;
\item $P_{U_i/W} (- T_i)$ for all $i$;
\item $J_{U_i/W} \cdot T_i$, for all $i$.
\end{enumerate}
\end{lem}
\proof
When $m = 1$, this is Theorem 1, Appendix of \cite{keel}.
Assume the statement is true for $m=r$. Consider the blowup
$$\Bl_{\{U_i: 1 \le i \le r+1 \}} W \rightarrow \Bl_{\{U_i: 1 \le i \le r \}} W $$
along the proper transform $U^r_{r+1}$ of $U_{r+1}$. Observe that
 $$P_{U^r_{r+1}/\Bl_{\{U_i: 1 \le i \le r \}} W} = P_{U_{r+1}/W}$$ and
$$J_{U^r_{r+1}/\Bl_{\{U_i: 1 \le i \le r \}} W} = (J_{U_{r+1}/W}, U^{r+1}_1, \ldots, U^{r+1}_r)$$
where $U^{r+1}_1, \ldots, U^{r+1}_r$ are the exceptional divisors of
$$\Bl_{\{U_i: 1 \le i \le r+1 \}} W \rightarrow W$$ corresponding to $U_1, \ldots, U_r$.
Hence the case of $r+1$ follows from the inductive assumption and Theorem 1 (i.e., $m=1$), Appendix of \cite{keel}.
\endproof

We identify $A^\bullet (W)$ as a subring of $A^\bullet (\Bl W)$ by means of
the injection $p^*: A^\bullet (W) \rightarrow A^\bullet (\Bl W)$ where
$p$ is the  projection $\Bl W \rightarrow W$.

\begin{lem} (\cite{FM})
Assume that $U$ and $V$ are smooth closed  subvarieties of $W$ and meet cleanly in a smooth
closed subvariety $Z$.
\label{Plist}
\begin{enumerate}
\item $P_{\Bl_Z U / \Bl_V W} (t) = P_{U/W} (t)$;
\item $P_{\Bl_Z U /\Bl_Z W} (t) = P_{U/W} (t - \widetilde{Z})$
where $\widetilde{Z}$ is the exceptional divisor
in $\Bl_Z W$;
\end{enumerate}
\end{lem}
\proof
This is basically  Lemma 6.2 of \cite{FM} except that
$U$ and $V$ meet cleanly instead of transversally.

(1) follows from that $N_{\Bl_Z U / \Bl_V W}$ is the pull-back of
$N_{U/W}$.

(2) follows from that $N_{\Bl_Z U / \Bl_V W} = p^* (N_{U/W}) \otimes \cO (- \widetilde{Z})|_{\Bl_Z U}$
and the verification used in \cite{FM}, where $p$ is the restriction to $\Bl_Z U$
of the map $\Bl_Z W \rightarrow W$.
\endproof

\begin{lem} (\cite{FM})
\label{Jlist}
Assume that $U$ and $V$ are smooth closed  subvarieties of $W$ and meet cleanly in a smooth
closed subvariety $Z$.
Assume also that $Z \hookrightarrow U, V \hookrightarrow W$ are all Lefschetz embeddings.
Then all the relevant inclusions below are Lefschetz embeddings, and
\begin{enumerate}
\item $J_{\Bl_Z U / \Bl_V W}  = J_{U/W} $ if $Z \ne \emptyset$;
\item $J_{\Bl_Z U / \Bl_V W}  = (J_{U/W}, \tilde{V})$ if $Z = \emptyset$, where $\tilde{V}$ is
the exceptional divisor in $\Bl_V W$;
\item $J_{\Bl_Z U /\Bl_Z W} = (J_{U/W},  [\Bl_Z V])$ if $Z \ne \emptyset$.
 Note that $ \Bl_Z V$ is the proper transform of $V$.
\end{enumerate}
\end{lem}
\proof
(1) and (2) together is Lemma 6.4 of \cite{FM}.

(3). By Lemma 6.5 of \cite{FM}, $J_{\Bl_Z U /\Bl_Z W} = (J_{U/W}, P_{V/W} ( - \widetilde{Z})$.
By Lemma \ref{Plist} (2), $P_{V/W} ( - \widetilde{Z}) = P_{\Bl_Z V /\Bl_Z W} (0) = [\Bl_Z V]$.
\endproof

\smallskip\noindent
{\sl Proof of Theorem \ref{ring}.}

\proof
First, fix some notation. Let $\cD_{=r}:=\{D_{r, 1} \ldots, D_{r, l_r} \}$ be the subset of
rank $r$ elements of the arrangement $\cD$.

We now prove the corresponding statement for $\Bl_{\cD_{\le r}} X$ by using induction on $r$.

When $r=0$, it follows directly from Lemma \ref{keel}.

Assume that the statement is true for $\Bl_{\cD_{\le r}} X$. That is,
the Chow ring $A^\bullet (\Bl_{\cD_{\le r}} X)$ is isomorphic to the  polynomial ring
$$A^\bullet (X) [T_{1, l_1}, \ldots T_{r, l_r}]/I_r$$
where $T_{m, j}$ corresponds to
{\it the exceptional divisor}  $E^{m+1}_{m, j}$ (see Theorem \ref{divisors} (1) for
the description of $E^{m+1}_{m, j}$) for
$D_{m,j} \in \cD_{\le r}$ and $I_r$ is the ideal generated by
\begin{enumerate}
\item $T_{m,i} \cdot T_{n,j}$, for incomparable $D_{m,i}$ and $D_{n,j}$ where $m, n \le r$;
\item $P_{D_{m, i}/X} (- \sum_{D_{n,j} \le D_{m,i}} T_{n,j})$ for all $(m, i)$ where $m \le r$;
\item $J_{D_{m, i}/X} \cdot T_{m, i}$; for all $i$ where $m \le r$.
\end{enumerate}

(Each variable $T_i$ geometrically corresponds to the cycle of the explicit
exceptional divisor $E^{r+1}_i$
in the blowup stage $\Bl_{\cD_{\le r}} X \rightarrow X$ where
$r= \hbox{rank} (D_i)$. This has to be beared in mind in the course of the rest of
the proof.)

Consider now the blowup
 $\Bl_{\cD^{r+1}_{\le 0}} (\Bl_{\cD_{\le r}} X) \rightarrow \Bl_{\cD_{\le r}} X$.
By Lemma \ref{keel}, $A^\bullet (\Bl_{\cD_{\le r+1}} X)$ is isomorphic to the polynomial ring
$$A^\bullet (\Bl_{\cD_{\le r}} X) [T_{{r+1},1}, \ldots, T_{r+1,l_{r+1}}]/I_{r+1}'$$
where $T_{{r+1},i}$ corresponds to the
exceptional divisor $E^{r+2}_{r+1,i}$ for $D_{r+1,i}$ of rank $r+1$ and $I_{r+1}'$ is generated by
\begin{enumerate}
\item $T_{r+1,i} \cdot T_{r+1,j}$, for $i \ne j$;
\item $P_{D^{r+1}_{r+1,i}/\Bl_{\cD_{\le r}} X} (- T_{r+1,i})$ for all $i$;
\item $J_{D^{r+1}_{r+1,i}/\Bl_{\cD_{\le r}} X} \cdot T_{r+1,i}$, for all $i$.
\end{enumerate}

For relation (2), we have $$P_{D^{r+1}_{r+1, i}/\Bl_{\cD_{\le r}} X} (- T_{r+1, i})
 = P_{D_{r+1, i}/X} (- \sum_{D_{m,j} \le D_{r+1, i}} T_{m,j}).$$
This is because by using Lemma \ref{Plist} (1) and (2)
repeatedly
 $$P_{D^{r+1}_{r+1, i}/\Bl_{\cD_{\le r}} X} (- T_{r+1, i}) =$$
 $$P_{D^r_{r+1, i}/\Bl_{\cD_{\le r-1}} X} (- T_{r+1, i} -
 \sum_{D_{r,j} < D_{r+1,i}} T_{r,j}) $$
$$= \ldots = P_{D_{r+1, i}/X} (- \sum_{D_{m,j} \le D_{r+1, i}} T_{m,j}).$$

For (3), we have by Lemma \ref{Jlist} (2) and (3) that
 $J_{D^{r+1}_{r+1,i}/\Bl_{\cD_{\le r}} X}$ is generated by
$$J_{D^r_{r+1, i}/\Bl_{\cD_{\le r-1}}} X, T_{r, j}, D^{r+1}_{r+1,l}$$
where $D_{r,j}$ and $D_{r+1, i}$ are incomparable and $D_{r+1, i} \cap D_{r+1, l}
= D_{r, h}$ for some $h$. Here $T_{r,j}$ comes from Lemma \ref{Jlist} (2),
while the proper transform $D^{r+1}_{r+1,l}$ presents due to  Lemma \ref{Jlist} (3).
But by the projection formula, the relation $D^{r+1}_{r+1,l} \cdot T_{r+1, i}$ follows
from Relation (1)  $T_{r+1,l} \cdot T_{r+1,i}$ and is thus redundant.
Hence, $I_{r+1}'$ is generated by
\begin{enumerate}
\item $T_{r+1,i} \cdot T_{m,j}$,  for incomparable $D_{r+1, i}$ and $D_{m,j}$, $ r \le m \le r+1$;
\item $P_{D_{r+1, i}/X} (- \sum_{D_{m,j} \le D_{r+1, i}} T_{m,j})$ for all $i$;
\item $J_{D^r_{r+1, i}/\Bl_{\cD_{\le r-1}} X} \cdot T_{r+1,i}$, for all $i$.
\end{enumerate}

The same argument as above used again shows that
 $J_{D^r_{r+1, i}/\Bl_{\cD_{\le r-1}} X}$ is generated by
 $$J_{D^{r-1}_{r+1, i}/\Bl_{\cD_{\le r-2}}}, T_{r-1, j}, D^r_{r,l}$$
where $D_{r-1,j}$ and $D_{r+1, i}$ are incomparable and $D_{r+1, i} \cap D_{r, l}
= D_{r-1, h}$ for some $h$.  Note that  $D_{r+1, i}$ and
$D_{r, l}$ are necessarily incomparable.
Again, the relation $D^r_{r,l} \cdot T_{r+1, i}$ follows
from $T_{r,l} \cdot T_{r+1, i}$ by the projection formula and is therefore redundant.
This reduces the relations above to
\begin{enumerate}
\item $T_{r+1,i} \cdot T_{m,j}$,
 for incomparable $D_{r+1, i}$ and $D_{m,j}$,  $r-1 \le m \le r+1$;
\item $P_{D_{r+1, i}/X} (- \sum_{D_{m,j} \le D_{r+1, i}} T_{m,j})$ for all $i$;
\item $J_{D^{r-1}_{r+1, i}/\Bl_{\cD_{\le r-2}} X} \cdot T_{r+1,i}$, for all $i$.
\end{enumerate}
 Keep repeating this procedure, we will eventually achieve
that
 $I_{r+1}'$ is generated by
\begin{enumerate}
\item $T_{r+1,i} \cdot T_{m,j}$,  for incomparable $D_{r+1, i}$ and $D_{m,j}$
      where $m \le r+1$;
\item $P_{D_{r+1, i}/X} (- \sum_{D_{m,j} \le D_{r+1, i}} T_{m,j})$ for all $i$;
\item $J_{D_{r+1, i}/ X} \cdot T_{r+1,i}$, for all $i$.
\end{enumerate}

Combine this with the inductive assumption on the case of $r$, the case of $r+1$
follows.

The statement of the theorem is the case when $r=$ rank of $\cD$.
\endproof

The same results, with essentially the same proof, hold for $H^* (\Bl_{\cD} X)$.

As pointed out in \cite{FM},
once an open variety $X^0$ is campactified by adding normal crossing divisors,
there is a standard way to construct a differential graded algebra $(\cA^\bullet, \diff)$.
A general theorem of Morgan asserts that $(\cA^\bullet, \diff)$ is a model for
the space $X^0$. Thus $(\cA^\bullet, \diff)$ determines the rational homotopy type
of $X^0$. In particular, $H^* (\cA^\bullet) = H^* (X^0)$.
This provides a practicable approach to compute, for example,  the cohomology
of the complements of arrangements in ${\Bbb C}^n$ or ${\Bbb P}^n$.

\section{Hodge polynomial of $\Bl_\cD X$}
\label{sectionhodge}

For any quasi projective variety $V$ there is a (virtual Hodge) polynomial $\e (V)$ in two variables
$u$ and $v$ which is uniquely determined by the following properties
\begin{enumerate}
\item If $V$ is smooth and projective, then $\e (V) = \sum h^{p, q} (-u)^p (-u)^q$;
\item If $U$ is a closed subvariety of $V$, then $\e (V) = \e (V     \setminus  U) + \e (U)$;
\item If $V \rightarrow B$ is a Zariski locally trivial bundle with fiber $F$, then
$\e (V) = \e (B) \e (F)$.
\end{enumerate}

The virtual Poincar\'e polynomial $\Po (V)$ is defined similarly.


\smallskip\noindent
{\sl Proof of Theorem \ref{hodge}.}

\proof We use induction on the number of elements in the arrangement of the subvarieties.

When $|\cD| = 1$, it is standard. Assume that the formula is true when
the number of elements in the arrangement of the subvarietie is less than $|\cD|$.

From the blowup
$$\Bl_\cD X = \Bl_{\cD^k_{\le 0}} (\Bl_{\cD_{\le k-1}} X) \rightarrow \Bl_{\cD_{\le k-1}} X$$
and the descriptions of the proper transforms of $D_i$ (i.e., Theorem \ref{transforms}),
we have
$$\e (\Bl_\cD X) = \e (\Bl_{\cD_{\le k-1}} X) + \sum_{\hbox{rank} (D_i) = k} \e
(\Bl_{D_{< D_i}} D_i) {{(uv)^{\dim X - \dim D_i} - uv} \over {uv -1}}.$$
By the inductive assumption,
$$\e (\Bl_{\cD_{\le k-1}} X) = \e (X) + \sum_{{\scriptsize \begin{array}{c}
D_{i_1} < \ldots <D_{i_{r+1}} \cr
 D_{i_{r+1}} := X \end{array}}}
\e (D_{i_1}) \prod_{j=1}^r {{(uv)^{\dim D_{i_{j+1}} - \dim D_{i_j}}} - uv \over {uv -1}}$$
where each $D_{i_j}$ is of rank $\le k-1$, and
$$\e (\Bl_{D_{< D_i}} D_i)= \e (D_i) + \sum_{{\scriptsize \begin{array}{c}
D_{i_1} < \ldots  < D_{i_r} \cr
D_{i_r}:= D_i \end{array}}}
\e (D_{i_1}) \prod_{j=1}^{r-1} {{(uv)^{\dim D_{i_{j+1}} - \dim D_{i_j}}} - uv \over {uv -1}}.$$
The formula in Theorem \ref{hodge} then follows from a direct computation from here.
Note that the convention $D_{i_{r+1}} :=X$ in the index of summation is a
manipulation to make the formula uniform and concise.

When $uv$ is substituted by $t^2$, the essentially same proof yields the formula for Poincar\'e
polynomials.
\endproof

\medskip\noindent
{\sl Proofs of Corollaries \ref{config} --  \ref{inCn}.}

\proof
Corollary \ref{config} is immediate from Theorem \ref{hodge}.

Corollary \ref{hyperplane} follows from Theorem \ref{hodge}
but needs a little manipulation of indexes to absorb the extra term $\e ({\Bbb P}^n)$
into the summation as indexed by ``$\emptyset < X$''. The index ``$\emptyset < H_{i_1}$'' is used
the same way for the factor $\e (H_{i_1})$.

Corollary \ref{m0nhodge} follows directly from \ref{hyperplane}.
The convention $|S_{i_{r+1}}| := n+1$ in the index is for a unified look of
the factors in the product.   (Note that
a {\it recursive} formula for
the Betti numbers of $\overline{M}_{0,n}$ was calculated by Keel via a different sequence of
blowups \cite{keel}.)

Corollary \ref{inCn}  is a special case of Theorem \ref{hodge}. Similar manipulation
of indexes as for Corollary \ref{hyperplane} is used.
\endproof

\section{Space of maps ${\Bbb P}^1 \rightarrow {\Bbb P}^n$}
\label{section:holomaps}

From the introduction,  $N_d ({\Bbb P}^n)$
is the space of equivalence classes of $(n+1)$-tuples $(f_0, \ldots, f_n) \ne 0$
where $f_i$ are homogeneous polynomials of degree $d$ in two variables, and
$$(f_0, \ldots, f_n) \sim (f'_0, \ldots, f'_n)$$
if $(f'_0, \ldots, f'_n)=c (f_0, \ldots, f_n)$ for some constant $c \ne 0$.

Let ${\bf 1}_d$ be the top partition $1+1+ \ldots +1$ of $d$ and
${\bf 0}_d$ the bottom (non) partition $d$ of $d$.
 If $\tau$ and $\tau'$ are
partitions of $d$, let $\tau \vee \tau'$ be the least partition that both $\tau$ and $\tau'$ proceed.
$\tau ({\bf 1}_r)$ denotes the partition $\tau + {\bf 1}_r$ of $d +r$.
$\rho (\tau)$ denotes the number of integers in the partition $\tau$.

Let $W^{(d)}$ be the $d$-th symmetric product of a variety $W$.
Then we should have
\begin{enumerate}
\item $N_{d_0, \tau} \cap  N_{ d_0, \tau'} = N_{d_0, \tau \vee \tau'}$;
\item $N_{d_0, \tau} \cap N_{ d'_0, \tau'} = N_{d_0', \tau ({\bf 1}_{d_0 - d_0'}) \vee \tau'}$
 if $d_0 > d_0'$.
\end{enumerate}

(1) can be reduced to configuration of
 unordered $(d -d_0)$ points in ${\Bbb P}^1$.
There, the result is clear.

(2) follows from the observation that
by factoring $d_0-d_0'$ linear factors from degree $d_0$ polynomials
we obtain $$N_{d_0, \tau} \cap N_{ d'_0, \tau'} = N_{d_0', \tau({\bf 1}_{d_0 - d_0'})}
\cap N_{ d'_0, \tau'}.$$ Hence it is reduced to (2).

As a poset, one can check the following partial relation:
$N_{d_0, \tau} > N_{ d'_0, \tau'}$ if one of the following holds
\begin{enumerate}
\item $d_0 = d_0'$, $\tau > \tau'$;
\item $d_0 > d_0'$, $\tau ({\bf 1}_{d_0 - d_0'}) \ge \tau'$.
\end{enumerate}

Unfortunately, it can be checked that $N_{d_0, \tau}$ may in general have singularities along lower strata.

The poset $\cN$ has the smallest element $N_{0, {\bf 0}_d}$. 
This is a smooth subvariety in $N_d ({\Bbb P}^n)$.
The strata of rank 2 have singularities along this subvariety, in general. However, 
it is possible that after blowing up
$N_{0, {\bf 0}_d}$, the singularites  the strata of rank 2 get resolved and their proper
transforms become separated so that the blowups along these proper transforms can be carried out.
Of course, the proper transforms of the strata on level 3 may still have singularities along
the proper transforms of the strata of level 2, the hope is that after blowing up
the proper transform of the strata of level 2,
they too get resolved and become separated so that the process can be carried on and on.

This requires an intensive analysis of singularities of the strata and the effects of blowups on them.
But the problem seems very interesting and of independent value, and calls for immediate investigation.






If the above turns out true,  then one can still compute the Hodge numbers.

For a partition $\tau$, let $\rho (\tau)$ be the number of integers in the partition. Then the same method
applied earlier would give that 
the polynomials
$\e (\Bl_\cN N_d ({\Bbb P}^n))$ and $\Po (\Bl_\cN N_d ({\Bbb P}^n))$ are given respectively by
$$\sum_{{\scriptsize \begin{array}{c}
 (d_{i_0}, \tau_{i_0}), (d_{i_0}', \tau_{i_0}') < (d_{i_1}, \tau_{i_1}) \ldots < (d_{i_{r+1}}, \tau_{i_{r=1}}) \cr
 d_{i_0}:=-1, \tau_{i_0} := \tau_{i_1} \cr
 d_{i_0}':=d_{i_1},  \tau_{i_0}' := \emptyset \cr
 d_{i_{r+1}}:=d, \tau_{i_{r+1}} := \emptyset  \end{array}}}
 \prod_{j=0}^r {{(uv)^{ (n+1)(d_{i_{j+1}} - d_{i_j}) + \rho (\tau_{i_{j+1}}) - \rho (\tau_{i_j})} - uv} \over {uv -1}},$$
$$
\sum_{{\scriptsize \begin{array}{c}
 (d_{i_0}, \tau_{i_0}), (d_{i_0}', \tau_{i_0}') < (d_{i_1}, \tau_{i_1}) \ldots < (d_{i_{r+1}}, \tau_{i_{r=1}}) \cr
 d_{i_0}:=-1, \tau_{i_0} := \tau_{i_1} \cr
 d_{i_0}':=d_{i_1},  \tau_{i_0}' := \emptyset \cr
 d_{i_{r+1}}:=d, \tau_{i_{r+1}} := \emptyset  \end{array}}}
 \prod_{j=0}^r {{(t^2)^{ (n+1)(d_{i_{j+1}} - d_{i_j}) + \rho (\tau_{i_{j+1}}) - \rho (\tau_{i_j})} - t^2} \over {t^2 -1}},$$
where $\rho(\emptyset) := 0$.

If $X$ is embedded in ${\Bbb P}^n$ such that the closure $N_d(X)$ in $N_d ({\Bbb P}^n)$
of the space  $N_d^0 (X)$ of
holomorphic maps of degree $d$ from ${\Bbb P}^1$ to $X$  is orbifold
and  meets nicely with
the subvarieties $N_{d_0, \tau}$
of the above arrangement $\cN$, then the blowup $\Bl_{\cN} N_d ({\Bbb P}^n)$ might induce
a blowup of $N_d(X)$ and
a nice  projective compactification $\Bl_{\cN (X)} N_d (X)$ of
 $N_d^0 (X)$ by adding normal crossing divisors, where $\cN (X)$ is the arrangement of the
subvarieties that are intersections of $N_d(X)$ and $N_{d_0, \tau}$.
We wonder to what extent this is the case for homogeneous spaces
 or more  generally convex varieties.

Replacing the two variables ($w_0, w_1$) by multiple variables, a formal
extension of the results
to maps ${\Bbb P}^m \rightarrow {\Bbb P}^n$ may be possible.

\section{Partial desingularization of GIT quotients}
\label{GIT}

In this section, the base field is assume to be of characteristics 0.

Let a reductive algebraic group $G$ act algebraically on a smooth projective variety $X$.
Let $L$ be a linearized ample line bundle over $X$ and $X^s = X^s (L)$ ($X^{ss} = X^{ss} (L)$)
the open subset of (semi) stable points in $X$.
 $X^{ss}     \setminus  X^s$ may not be empty.

Replacing $L$ by a large tensor power we may assume that $L$ is very ample and hence
induces an equivariant embedding $X \hookrightarrow {\Bbb P}^n$.

Let $\Re$ be the set of conjugacy classes of all connected reductive subgroups of $G$ and
 $R$  a representative of an arbitrary class in $\Re$. Define
$$Z_R^{ss} = \{ [x_0, \ldots, x_n] \in X | (x_0, \ldots, x_n) \; \hbox{is fixed by} \; R \}\cap X^{ss} .$$
Then $GZ_R^{ss}$ is a  closed smooth subvariety in $X^{ss}$ by 5.10 and 5.11 of \cite{Kirwan}.
Different $Z_R^{ss}$, as connected components of reductive subgroups, meet cleanly. So do the corresponding
$GZ_R^{ss}$ by 5.10 of \cite{Kirwan}.

\begin{lem} Let $\cR$ be the set $\{GZ_R^{ss} : R \in \Re \}$. Then $\cR$ is an arrangement
 of smooth subvarieties of $X^{ss}$.
\end{lem}

 Note that the above statement is void when $X^{ss} = X^s$.

\begin{lem}
\label{quotientnc} Let $W$ be a smooth  algebraic variety acted on
by a  reductive algebraic group $G$. Let $W^0$ be a $G$-invariant open
subset such that  $W     \setminus  W^0 = \bigcup_i D_i$ is
a union of smooth $G$-invariant divisors $D_i$  with normal crossings.
Assume that $L$ is a linearization of the $G$-action such that
$W^{ss} = W^s$. Then $W^{s}/G   \;\;  \setminus  \;\; (W^{s} \cap W^0)/G = \bigcup_i  (W^{s} \cap D_i)/G$
is a union of
normal crossing divisors with at worst finite quotient singularities\footnote{
Such divisors are said to meet transversally if up to a finite etal\'e covering,
they meet transversally in the usual sense. }.
 (Note that some of $(W^{s} \cap D_i)/G$ may be empty.)
\end{lem}

\proof
It suffices to check that $\{(W^s \cap D_i)/G \}_i$ meet transversally. This is a
local question.

Given a point $x \in W^s$. By Luna's \'etale slice theorem, there is a
$G$-invariant open neighborhood $W_x$ of the point $x$ in $W^s$ and a
smooth  $G_x$-invariant subvariety $S_x$ in $W_x$ containing $x$ such that $W_x = G \cdot S_x$ and
the natural map
$$  G \times_{G_x} S_x \rightarrow G \cdot S_x = W_x$$ is \'etale.
Here $G_x$ is the finite isotropy subgroup of $G$ at $x$.
This induces \'etale maps
$$ G \times_{G_x}(S_x \cap D_i)  \rightarrow G \cdot (S_x \cap D_i) = (G \cdot S_x) \cap D_i =W_x \cap D_i $$
where $S_x \cap D_i$ is either empty or a divisor in $S_x$.
Now, since $\{D_i\}$ meet transversally in $W_x$, we have that their corresponding
$\{D_i \cap S_x\}$ meet transversally in $S_x$. 
Hence the quotients $\{(S_x \cap D_i)/G_x\}$ by the finite group $G_x$, 
as divisors with at worst finite quotient singularities,
 meet transversally by definition (see the footnote 3).
Finally, we can use the natural identification $(W_x \cap D_i)/G \cong  (S_x \cap D_i)/G_x$ 
(local analytically)
to conclude the proof.
\endproof

 Let $\cD$ be any arrangement of $G$-invariant subvarieties in $X^{ss}$
such  that $\cD$
contains $\cR$ as a subarrangement.
Let $E$ be the exceptional divisors of $p: \Bl_\cD X^{ss} \rightarrow X^{ss}$.
$M_d = p^*(L^{\otimes d}) \otimes \cO(- E)$ admits a linearization and
$(\Bl_\cD X^{ss})^{ss}(M_d)$ is independent of  sufficiently large $d$. In the following,
the stability and quotient are taken with respect to the linearization $M_d$ for a
fixed sufficiently large $d$.
By Lemma 6.1 of \cite{Kirwan}, $(\Bl_\cR X^{ss})^{ss} = (\Bl_\cR X^{ss})^s$.
Then by  relative GIT (e.g, Theorems 3.11 and 4.4 of \cite{Hu96}),
 $(\Bl_\cD X^{ss})^{ss} = (\Bl_\cD X^{ss})^s$.


\begin{cor}
\label{useform0nbar} The variety $\Bl_\cD X^{ss}$ has a  geometric quotient such that
the following diagram is commutative
\begin{equation*}
\begin{CD}
 (\Bl_\cD X^{ss})^s  @>{\hat{q}}>> (\Bl_\cD X^{ss})^s /\!/ G\\
@V{p}VV @V{p_G} VV \\
 X^{ss}  @>q>> X^{ss}/\!/G
\end{CD}
\end{equation*}
 and the complement of
$(\hat{q} \circ p^{-1} (X^s \cap X^0))/\!/G \cong (X^s \cap X^0)/\!/G$ in $(\Bl_\cD X^{ss})^s /\!/ G $
is a union of normal crossing divisors with at worst finite quotient singularities,
 where $X^0 = X^{ss}     \setminus  \bigcup \cD$.
Moreover, $$p_G: (\Bl_\cD X^{ss})^s /\!/ G \rightarrow  X^{ss}/\!/G$$
is a blowup along the induced arrangement of the images of the subvarieties in $\cD$.
\end{cor}
\proof
This follows from the combination of Theorem \ref{mainlemma} and Proposition 6.9 of \cite{Kirwan}.
One needs to observe that each stage of the blowups $p: \Bl_\cD X^{ss} \rightarrow X^{ss}$
yields (by applying Lemma 3.11 of \cite{Kirwan})
a corresponding stage of blowups $p_G: (\Bl_\cD X^{ss})^s /\!/ G \rightarrow  X^{ss}/\!/G$
and  once $D_i$ and $D_j$ get separated in  certain stage the same become true for their
images in $X^{ss}/\!/G$.
\endproof

Note the blowup $p_G$ is not covered by Theorem \ref{mainlemma} due to the presence
of singularities. But the blowing up procedure and the reason that it can be carried out
is essentially the same,  as indicated in the above proof.

When $\cD = \cR$, we recover Kirwan's partial desingularization of $ X^{ss}/\!/G$ (\cite{Kirwan}).

\medskip

\bibliographystyle{amsplain}
\makeatletter \renewcommand{\@biblabel}[1]{\hfill#1.}\makeatother

\vskip .4cm
\noindent
Department of Mathematics, Unversity of Arizona, Tucson, AZ85721
yhu@@math.arizona.edu

\end{document}